\documentclass[12pt,letterpaper]{article}
\usepackage{lipsum}
\usepackage[margin=1.5in]{geometry}
\usepackage{tikz}
\usetikzlibrary{arrows.meta, positioning}
\usepackage[applemac]{inputenc}
\usepackage{amsmath,amssymb}
\usepackage{amsfonts}
\usepackage{amsmath,amssymb}
\usepackage[applemac]{inputenc}
\usepackage{color}
\usepackage{color, colortbl}
\usepackage{amsmath}
\usepackage{amssymb}
\usepackage{bbm}
\usepackage{graphicx}
\usepackage{tikz}
\usepackage{geometry}
\usetikzlibrary{arrows}
\usepackage{amsfonts}
\usepackage{amsmath,amssymb}
\usepackage[applemac]{inputenc}
\usepackage{color}
\usepackage{amsmath}
\usepackage{amssymb}
\usepackage{graphicx}
\usepackage{tikz}
\newtheorem{theorem}{Theorem}[section]
\newtheorem{lemma}[theorem]{Lemma}
\newcommand{\E}{E}

\newtheorem{corollary}[theorem]{Corollary}
\newtheorem{definition}[theorem]{Definition}

\newtheorem{example}[theorem]{Example}

\newtheorem{problem}[theorem]{Problem}
\newtheorem{remark}[theorem]{Remark}

\usetikzlibrary{arrows,shapes,automata,petri}

\newcommand{\dproof}{\noindent {Proof.} \quad}
\newcommand{\fproof}{\hfill $\square$ \bigskip}

\numberwithin{equation}{section}

\definecolor{LightCyan}{rgb}{0.88,1,1}

\def\RR{{\mathbb{ R}}} 

\def\EE{{\mathbb{ E}}}

\def\1B{\text{1\!\!I}}

\def\<{\langle}
\def\>{\rangle}

\def\E{\mathbb{E}}
\def\R{\mathbb{R}}

\begin{document}
\title{Fokker-Planck equations for conditional McKean-Vlasov systems driven by Brownian sheets}
\author{N. Agram$^{1},$ B. \O ksendal$^{2}$, F. Proske$^{2}$ and O. Tymoshenko$^{2,3}$}
\date{\today}
\maketitle

\footnotetext[1]{Department of Mathematics, KTH Royal Institute of Technology 100 44, Stockholm, Sweden. 
Email: nacira@kth.se. Work supported by the Swedish Research Council grant (2020-04697),the
Slovenian Research and Innovation Agency, research core funding No.P1-0448.}

\footnotetext[2]{%
Department of Mathematics, University of Oslo, Norway. \\
Emails: oksendal@math.uio.no, proske@math.uio.no, otymoshenkokpi@gmail.com}
\footnotetext[3]
{Department of Mathematical Analysis and Probability Theory, NTUU Igor Sikorsky Kyiv Polytechnic Institute, Kyiv, Ukraine.  MSC4Ukraine grant (AvH ID:1233636).}

\begin{abstract}
We investigate conditional McKean-Vlasov equations driven by time-space white noise, motivated by the propagation of chaos in an \(N\)-particle system with space-time Ornstein-Uhlenbeck dynamics. The framework builds on the stochastic calculus of time-space white noise, utilizing tools such as the two-parameter It\^o formula, Malliavin calculus, and orthogonal decompositions to analyze convergence and stochastic properties. Existence and uniqueness of solutions for the associated stochastic partial differential equations (SPDEs) are rigorously established. Additionally, an integral stochastic Fokker-Planck equation is derived for the conditional law, employing Fourier transform methods and stochastic analysis in the plane. The framework is further applied to a partial observation control problem, showcasing its potential for analyzing stochastic systems with conditional dynamics.
\end{abstract}
\paragraph{MSC(2020):}

60G60; 60G35; 60H15; 60H20; 60H40; 91A15; 93E11; 93E20.\\

\textbf{Keywords:} McKean-Vlasov (mean-field) SPDE, Brownian sheet, conditional probability distribution, stochastic Fokker-Planck equation, partial observation control. 
\section{Introduction}
To the present time, various studies have been published regarding stochastic integrals in the plane,  (see, for
example, \cite{CW}, \cite{I2},\cite{MiS}, \cite{WZ74}). 
 Furthermore, the studies referenced in \cite{NualartSanz}, \cite{FarreNualart},  and \cite{Yeh81}  provide a foundation for understanding the conditions under which both strong and weak solutions exist for stochastic differential equations (SDEs) on the plane, alongside requirements for smoothness and other properties.
However, significant unresolved issues remain in the theory of stochastic differential equations on the plane. Additionally, McKean-Vlasov stochastic partial differential equations (SPDEs), which arise as the mean-field limit of stochastic partial differential equation systems, have gained increasing attention in recent years because of their potential applications in areas such as,  for example, neurophysiology, quantum field theory, kinetic theory, finance and other economics fields. It is important to note that this type of equation in the unconditional case was first investigated by H. McKean in \cite{McK}. However, there is a lack of results in the existing literature regarding  mean-field type SPDEs driven by Brownian sheets.

The one-parameter case of McKean-Vlasov SDEs in infinite dimensions have been studied recently by Hong \textit{et al.}  \cite{HLL}. Weak and strong existence and weak and strong uniqueness results for
multi-dimensional stochastic McKean-Vlasov equations are established under linear growth or integrability conditions in \cite{MiSV}. As for existence and uniqueness
of the solution of McKean-Vlasov SDEs, and the associated Fokker-Planck equation,
we refer the reader to Bogachev \textit{et al.}  \cite{BkR} and Barbu \& R{\"o}ckner \cite{BR}. In the one-parameter case, Fokker-Planck equations in infinite dimensions have been considered in  \cite{BPR}, where the authors introduced a general method for establishing the uniqueness of solutions in infinite-dimensional spaces.  These equations have also been studied by  Agram \textit{et al.}  \cite{AO, APO, AR} with applications extending to optimal control and even deep learning. Recently, Agram \textit{et al.}  \cite{AOPT1} have considered various applications, including the optimal control of time-space SPDEs driven by a Brownian sheet, but no mean-field terms are included. In \cite{AOPT2} authors state and prove the Fokker-Planck equation and illustrate 
the results for some time-space SPDEs.

The purpose of this paper is to generalise the results in \cite{AOPT2}, and study the conditional law of the solution of McKean-Vlasov type SPDEs driven by Brownian sheets. Specifically, we analyze the evolution of conditional distributions at a specific time $t$ and space $x$, denoted by $\mu_{t,x}$   of solutions to McKean-Vlasov SPDE with common noise of the form
\small
\begin{align}\label{Y2}
Y(z)&=\left[\begin{array}{clcr}
Y_1(z)\\Y_1(z)\\...\\Y_n(z)
\end{array} \right]
=\left[\begin{array}{clcr}
Y_1(0) \\Y_2(0)\\ ... \\Y_n(0)
\end{array} \right] 
+\int_{R_z} \left[ \begin{array}{c}
\alpha_1(\zeta)\\ \alpha_2(\zeta)\\ ...  \\\alpha_n(\zeta)
\end{array} \right] d\zeta
+\int_{R_z}\left[ \begin{array}{rc}
\beta_1(\zeta) \\ \beta_2(\zeta) \\ ... \\ \beta_n(\zeta) 
\end{array} \right]
\left[\begin{array}{rc} B_1(d\zeta)\\ B_2(d\zeta) \\ ... \\B_m(d\zeta)
\end{array} \right] 
\end{align}
where $R_z=[0,t]\times[0,x]$, $z=(t,x)$,
$\alpha_k(z)=\alpha_k(z,Y(z),\mu_{z}) \in \mathbb{R} \text{ for all } k=1,2, ...,n,$ and 
$$\beta_{\ell}(z)=\beta_{\ell}(z,Y(z),\mu_{z})= (\beta_{\ell,1}(z), \beta_{\ell,2}(z), ... , \beta_{\ell,m}(z)) \in \mathbb{R}^{1 \times m}$$
 is a $m-$dimensional row vector for all $\ell = 1,2, ..., n $. 
Note that 
 $$B(z)=(B_1(z), B_2(z), ... ,B_m(z))^{T} \in \mathbb{R}^{m \times 1}; \quad z=(t,x)$$ is a $m$-dimensional Brownian sheet defined on a filtered probability space  $(\Omega, \mathcal{F} , \mathbb{F} = \{\mathcal{F}^{(m)}_{t,x}\}_{t,x\geq0}, \mathbb{P} )$. Assuming the same boundary conditions as in \eqref{Y2}, the corresponding equivalent differential form of \eqref{Y2} in terms of  time-space white noise  $ \overset{\bullet}{B}$ and Wick product $\diamond $ is%
\begin{equation}
\frac{\partial^{2}}{\partial t\partial x}Y(t,x)=\alpha(t, x,  Y
(t,x),\mu_{t,x})+\beta(t, x, Y(t,x),\mu_{t,x})\diamond\overset{\bullet}{B}(t,x)\label{(2)}%
\end{equation}
 with coefficients $\alpha=(\alpha_1, \alpha_2, ... ,\alpha_n)\,\, \textnormal{and}\, \beta=(\beta_{i,j})_{1\leq i \leq n, 1\leq j\leq m}$.
 
The identity of \eqref{Y2} and (\ref{(2)}) comes from the fact that%
\begin{equation*}
\int_{R_{(t,x)}}\varphi(s,a)B(ds,da)=\int_{R_{(t,x)}}\varphi(s,a
)\diamond\overset{\bullet}{B}(s,a)dsda,\text{ } \quad \forall\varphi, t, x.
\end{equation*}
See e.g. Holden \textit{et al.} \cite{HOUZ}.\\

The conditional probability distribution of $Y(t,x)$  at a given time $t$  and space  $x$  is represented as $\mu_{t,x}=\mathcal{L}(Y(t,x) |\mathcal{F}_{t,x}^{(1)}) $, where 
 $\mathcal{F}_{t,x}^{(1)}$ is the filtration generated by the first component $B_1,$  up to time 
$t$ and space $x$ (commonly referred to as common noise). The remaining components of the 
$m$-dimensional Brownian sheet $B$ correspond to idiosyncratic noises. This conditional distribution plays a crucial role in understanding the behavior of $Y(t,x)$  and its evolution, especially in the context of prediction and control under the influence of randomness based on that information.

In many real-world scenarios, systems consist of multiple interacting agents or components, influenced by shared sources of randomness or uncertainty. For instance, in financial and economic models, \textit{common noise} serves as a shared stochastic driver, influencing the behavior of assets, markets, or agents in a synchronized manner. Such modeling helps capture the interconnected dynamics of these agents, offering insights into their collective responses to external shocks. This framework aligns with the microeconomic theory of common shocks, enabling analysis of system-wide dynamics and decision-making in complex environments.

Moreover, the concept of common noise closely ties to filtering theory, where only partial observations of the system are available. This is evident when the filtration $\mathcal{F}_{t,x}^{(1)}$ is trivial, corresponding to cases without additional filtering, akin to non-common noise settings where all information is directly accessible. These cases have been examined in works like Agram \textit{et al.} \cite{AOPT2}, emphasizing scenarios where inference is unnecessary.

An applied example arises in environmental monitoring, such as tracking pollution levels in a city via a network of sensors. The true state of the system \( Y(z) \) evolves under a combination of deterministic and stochastic influences. To account for uncertainty stemming from incomplete information, measurement inaccuracies, or external disturbances, we incorporate the conditional measure \( \mu_z \), defined as \( \mu_z = \mu_{t,x} = \mathcal{L}(Y(t,x) | \mathcal{F}_{t,x}^{(1)}) \). This measure represents the conditional law of the state given the available noisy observations, capturing the probabilistic nature of the system.

The evolution of the state \( Y(z) \) is described by the following dynamics:
\[
Y(z) = Y(0) + \int_{R_z} \alpha(\zeta, Y(\zeta), \mu_\zeta, u(\zeta)) \, d\zeta + \int_{R_z} \beta(\zeta, Y(\zeta), \mu_\zeta, u(\zeta)) \, B(d\zeta),
\]
where \( \alpha \) and \( \beta \) represent the drift and diffusion terms, respectively, and \( u(z) \) denotes a control strategy. The uncertainty in the pollution level, captured by \( \mu_z \), arises from factors such as imperfect measurements, unmodeled external variables, or environmental fluctuations (e.g., weather or traffic).

The observation process provides noisy sensor data about the pollution level:
\[
dG(z) = dB_1(z),
\]
where \( G(z) \) is the noisy observation, and \( B_1(z) \) is an independent Brownian sheet. The task is to design a control strategy \( u(z) \), such as traffic rerouting or industrial emissions regulation, to minimize pollution levels while accounting for the uncertainty in the system state.

The control problem is formulated as a partial observation optimization problem:
\[
\sup_{u \in \mathcal{A}_G} J(u),
\]
where the performance functional \( J(u) \) is given by:
\[
J(u) = \mathbb{E}\left[\int_0^T \int_0^X \ell(\zeta, Y(\zeta), u(\zeta)) \, d\zeta + k(Y(T, X))\right].
\]
Here, \( \ell \) represents the running cost, \( k \) denotes the terminal cost, and \( \mathcal{A}_G \) is the set of admissible control strategies adapted to the filtration generated by the observation process \( G(z) \).

This formulation captures the essence of decision-making under uncertainty, where actions are based on indirect and noisy observations of the true state. The inclusion of \( \mu_z \) enables the model to account for the probabilistic evolution of the pollution level and facilitates the derivation of optimal strategies. By solving this partial observation control problem, the framework provides a systematic approach to addressing complex environmental management challenges where complete information is unavailable.

The paper is structured as follows:

In {Section 2} we review some preliminary concepts that will be used throughout this work. Specifically, we present some background about the stochastic calculus of time-space white noise, in particular, the two-parameter It\^o's formula. 

In {Section 3} to motivate conditional McKean-Vlasov SDEs in the two-parameter case, the propagation of chaos for space-time Ornstein-Uhlenbeck SDEs is analyzed by considering an $N$-particle system where interactions are described through averaged dynamics influenced by stochastic coefficients and Brownian motion. The equations are further studied using the Malliavin derivative and orthogonal decomposition to explore their stochastic properties.

In {Section 4} the existence and uniqueness results for the solution of the conditional McKean-Vlasov SPDE are proved.

In {Section 5} an integral stochastic Fokker-Planck equation is derived and proved for the conditional law of a time-space McKean-Vlasov equation, using techniques such as Fourier transforms of measures and It\^o's formula for stochastic integrals in the plane.

In {Section 6} we establish the corresponding SPDE for the conditional law.

In {Section 7} we provide an application to partial observation control.

\section{Background}
In this section we give some background about the associated stochastic calculus in the plane. It covers two-parameter stochastic integrals and their extensions for multi-dimensional processes. 
\subsection{Two-parameters stochastic calculus}
Throughout this work, we denote by $\{B(t,x); t,x \geq 0 \}$ a  Brownian sheet and $(\Omega, \mathcal{F}, P)$ a complete probability space on which we define the (completed) $\sigma$-field $\mathcal{F}_{t,x}$ generated by $ B(s,a), s \leq t, a \leq x$.
Wong \& Zakai \cite{WZ} generalized the notion of stochastic integrals with respect to one-parameter Brownian motion to stochastic intergrals driven by the two-parameter Brownian sheet. Let us denote by $\mathbb{R}^2_+$ the positive quadrant of the plane and let $z\in \mathbb{R}^2_+$.
In Cairoli \cite{Cairoli72} a first type stochastic integral with respect to the two-parameter Brownian motion is defined, denoted by:
\begin{align*}
   \int_{\mathbb{R}^2_+} \phi(z)B(dz),
\end{align*}
and a second type \cite{WZ74} stochastic integral denoted by
\begin{align*}
   \int_{\mathbb{R}^2_+} \int_{\mathbb{R}^2_+} \psi(z,z')B(dz)B(dz').
\end{align*}
These concepts provide the foundation for advanced stochastic calculus on the plane.

We recall a multi-dimensional It\^{o} formula for stochastic integrals in the plane:
\begin{theorem}
[Multi-dimensional It\^{o} formula, Wang \& Zakai \cite{WZ}]

Suppose
\small
\begin{align}\label{Y3}
Y(z)&=\left[\begin{array}{clcr}
Y_1(z)\\Y_1(z)\\...\\Y_n(z)
\end{array} \right]
=\left[\begin{array}{clcr}
Y_1(0) \\Y_2(0)\\ ... \\Y_n(0)
\end{array} \right] 
+\int_{R_z} \left[ \begin{array}{c}
\alpha_1(\zeta)\\ \alpha_2(\zeta)\\ ...  \\\alpha_n(\zeta)
\end{array} \right] d\zeta
+\int_{R_z}\left[ \begin{array}{rc}
\beta_1(\zeta) \\ \beta_2(\zeta) \\ ... \\ \beta_n(\zeta) 
\end{array} \right]
\left[\begin{array}{rc} B_1(d\zeta)\\ B_2(d\zeta) \\ ... \\B_m(d\zeta)
\end{array} \right],
\end{align}
where
\begin{itemize}
    \item $R_z=[0,t]\times[0,x]$, $z=(t,x)$, $\alpha_k(z) \in \mathbb{R} \text{ for all } k=1,2, ...,n;$  
    \item $\beta_{\ell}(z)= (\beta_{\ell,1}(z), \beta_{\ell,2}(z), ... , \beta_{\ell,m}(z)) \in \mathbb{R}^{1 \times m}$ 
 is an $m$-dimensional row vector for all  $\ell = 1,2, ..., n; $ 
 \item  $B(z)=(B_1(z), B_2(z), ... , B_m(z))^{T} \in \mathbb{R}^{m \times 1}$ is an $m$-dimensional Brownian sheet. 
\end{itemize}

Then, if $f:\mathbb{R}^n \mapsto \mathbb{C}$ is smooth, we have
\small
\begin{align*}\label{(17)}
&f(Y(z)) =f(Y(0))+\int_{R_{z}} \sum_{k=1}^n \frac{\partial f}{\partial y_k}(Y(\zeta))\Big[\alpha_k(\zeta
)d\zeta+\beta_k(\zeta)B(d\zeta)\Big]\\
&+\tfrac{1}{2}\int_{R_{z}}\sum_{k,\ell =1}^n \frac{\partial^2 f}{\partial y_k \partial y_{\ell}} (Y(\zeta))\beta_k(\zeta) \beta_{\ell}^{T}(\zeta)d\zeta\\
& +\iint\limits_{R_{z}\times R_{z}}\sum_{k,\ell =1}^n \frac{\partial ^2 f}{\partial y_k \partial y_{\ell}}(Y(\zeta\vee\zeta^{\prime
}))\beta_k(\zeta)B(d\zeta)\beta_{\ell}(\zeta^{\prime})B(d\zeta^{\prime})\\ 
&+\iint\limits_{R_{z}\times R_{z}}\Big\{\sum_{k,\ell =1}^n \frac{\partial ^2 f}{\partial y_k \partial y_{\ell}}(Y(\zeta\vee\zeta^{\prime
}))\beta_k(\zeta^{\prime}) \alpha
_{\ell}(\zeta)\\
&+\tfrac{1}{2}\sum_{k,\ell,p=1}^n \frac{\partial ^{(3)} f}{\partial y_k \partial y_{\ell} \partial y_p}(Y(\zeta\vee\zeta^{\prime}))
\beta_k(\zeta^{\prime}) \Big[\beta_{\ell}(\zeta)\beta_p^T (\zeta)\Big]\Big\}d\zeta B(d\zeta^{\prime})\\
&+\iint\limits_{R_{z}\times R_{z}}\Big\{\sum_{k,\ell =1}^n \frac{\partial ^2 f}{\partial y_k \partial _{\ell}}(Y(\zeta\vee\zeta^{\prime
}))\beta_k(\zeta)\alpha_{\ell}(\zeta^{\prime})\\
&+\tfrac{1}{2}\sum_{k,\ell,p=1}^n \frac{\partial ^{(3)} f}{\partial y_k \partial y_{\ell} \partial y_p}(Y(\zeta\vee\zeta^{\prime}))
\beta_k(\zeta) \Big[\beta_{\ell}(\zeta^{\prime})\beta_p^T (\zeta^{\prime})\Big]\Big\}B(d\zeta)d\zeta^{\prime}\\
& +\iint\limits_{R_{z}\times R_{z}}I(\zeta \bar{\wedge} \zeta^{\prime}) \Big\{\sum_{k,\ell =1}^n \frac{\partial ^2 f}{\partial y_k \partial y_{\ell}}(Y(\zeta\vee\zeta^{\prime
}))\alpha_k(\zeta^{\prime}) \alpha_{\ell}(\zeta)\\
&
+\tfrac{1}{2}\sum_{k,\ell,p=1}^n \frac{\partial ^{(3)} f}{\partial y_k \partial y_{\ell} \partial y_p}(Y(\zeta\vee\zeta^{\prime}))
\Big[\alpha_k(\zeta^{\prime}) \beta_{\ell}(\zeta)\beta_p^T (\zeta)+\alpha_k(\zeta) \beta_{\ell}(\zeta^{\prime}) \beta_p^{T}(\zeta^{\prime})\Big] \\
& +\tfrac{1}{4} \sum_{k,\ell,p,q=1}^n \frac{\partial ^4 f}{\partial y_k \partial y_{\ell} \partial y_p \partial y_q}(Y(\zeta\vee\zeta^{\prime}))\beta_k(\zeta^{\prime})\beta_\ell^{T}(\zeta^{\prime}%
)\beta_p(\zeta)\beta_q^{T}(\zeta)\Big\}d\zeta d\zeta^{\prime}.
\end{align*}
\end{theorem}



\section{Propagation of Chaos for space time Ornstein-Uhlenbeck SDEs}
  
In order to motivate conditional McKean-Vlasov SDEs in the two-parameter
case from the viewpoint of propagation of chaos, consider now the following
linear $N-$particle system, $i=1,...,N$, for $0\leq
t$, $x\leq T$,
\[
Y^{i,N}(t,x)=Y_{i}(0)+\int_{0}^{t}\int_{0}^{x}\left( \frac{1}{N}%
\sum_{j=1}^{N}a_{j}Y^{j,N}(s,a)-Y^{i,N}(s,a)\right) dsda+B_{i}(t,x),
\]%
where $a_{j}:\Omega \rightarrow \mathbb{R}$, $j=1,...,N$ are stochastic
coefficients and $B=(B_{1},...,B_{N})$ is a $N-$dimensional Brownian sheet.
We assume that $\Omega =\Omega _{1}\times \Omega _{2}$ for sample spaces $%
\Omega _{i}$, $i=1,2$ and $\mathbb{P}=\mathbb{P}_{1}\times \mathbb{P}_{2}$ for (complete) probability measures $\mathbb{P}_{i}$ on $\mathcal{F}%
_{i}$, $i=1,2$. Suppose that $a_{j}(\omega _{1},\omega
_{2})=b_{j}((W_{t,x}(\omega _{1}))_{0\leq t,x\leq T})$, $j=1,...,N$ and $%
B(\omega _{1},\omega _{2})=B(\omega _{2})$ for $\omega _{i}\in \Omega _{i}$, 
$i=1,2$, where $W$ is another $d-$dimensional Brownian sheet and $b_{j}$, $%
j=1,...,N$ are Borel measurable functions. So the coefficients $a_{j}$, $%
j=1,...,N$ are independent of $B$. Denote by $\mathcal{G}$ the sub-$\sigma -$%
algebra of $\mathcal{F}=\mathcal{F}_{1}\times \mathcal{F}_{2}$ generated by $%
W$.      

The system of equations above, which could, for example, describe the dynamics of interacting waves in the ocean, can be written as:
\[
Y(t,x)=Y(0)+\int_{0}^{t}\int_{0}^{x}\left( \frac{1}{N}A-I_{N}\right)
Y(s,a)dsda+B(t,x),
\]%
where $0\leq t,x\leq T$, $Y(0)\in \mathbb{R}^{N}$, $%
Y(t,x):=(Y^{1,N}(t,x),...,Y^{N,N}(t,x))$, $I_{N}$ is the unit matrix and%
\[
A:=\begin{pmatrix} a_{1} & \cdots & a_{N} \\ \vdots & & \vdots \\ a_{1} &
\cdots & a_{N} \end{pmatrix}.
\]%
By applying the Malliavin derivative $D_{u,v}$ in the direction of the
Brownian sheet $B$ (see \cite{Nualart},  \cite{DOP}) to both sides of the latter
equation, for $0\leq u\leq t$, $0\leq v\leq x$, we obtain that%
\[
D_{u,v}Y(t,x)=\int_{u}^{t}\int_{v}^{x}\left( \frac{1}{N}A-I_{N}\right)
D_{u,v}Y(s,a)dsda+\chi _{R(t,x)}(u,v)I_{N}
\]%
a.e. See (2.75) p. 144 in \cite{Nualart}. Using Picard iteration, we see that%
\begin{align*}
&D_{u,v}Y(t,x)=&\\&=\sum_{n\geq
0}\int_{R(u,v,t,x)}\int_{R(u,v,t_{1},x_{1})}...\int_{R(u,v,t_{n-1},x_{n-1})}%
\prod\limits_{j=1}^{n}\left( \frac{1}{N}A-I_{N}\right)
dt_{n}dx_{n}...dt_{1}dx_{1}\nonumber\\
&=\sum_{n\geq 0}\frac{1}{(n!)^{2}}\left( (t-u)(x-v)\left( \frac{1}{N}%
A-I_{N}\right) \right) ^{n}=f\left( (t-u)(x-v)\left( \frac{1}{N}%
A-I_{N}\right) \right) \text{,}
\end{align*}
a.e., where $R_{(u,v,t,x)}=\{(s,a);u\leq s\leq t,v\leq a\leq x\}$.\newline
Here $f$ is a function which is related to the Bessel function of order zero
and given by
$$\displaystyle f(y)=\sum_{n\geq 0}\frac{1}{(n!)^{2}}y^{n}\text{.}$$
On the other hand, we know that $Y(t,x)$ has $\mathbb{P}_{1}-$a.e. the
orthogonal decomposition 
\begin{equation*}
Y(t,x)=\mathbb{E}_{\mathbb{P}_{2}}\left[ Y(t,x)\right] +\sum_{n\geq
1}I_{n}(f_{n})  
\end{equation*}%
with respect to $B$ in $L^{2}(\left[ 0,T\right] ^{2}\times \mathbb{P}%
_{2})$ , where $f_{n}=f_{n}(\omega _{1})(\in L_{s}^{2}(\left[ 0,T\right]
^{2n})$ (the space of symmetric square integrable kernels), $(t,x)\in \left[
0,T\right] ^{2}$ and $I_{n}$ is the multiple Wiener-It\^{o} stochastic
integral with respect to the Brownian sheet $B$ (see \cite{Nualart}). Here $%
\mathbb{E}_{\mathbb{P}_{2}}$ denotes the expectation with respect to 
$\mathbb{P}_{2}$. Applying the Malliavin derivative, we find that 
$\mathbb{P}_{1}-$a.e. 
\[
D_{u,v}Y(t,x)=\sum_{n\geq 1}nI_{n-1}(f_{n}(z_{1},...,z_{n-1},z)
\]%
in $L^{2}(\left[ 0,T\right] ^{2}\times \mathbb{P}_{2})$, $z=(u,v)$. Then
taking the variance with respect to $\mathbb{P}_{2}$ on both sides of the
latter equation, entails that $\mathbb{P}_{1}-$a.e.%
\[
\sum_{n\geq 2}n^{2}(n-1)!\left\Vert f_{n}(z_{1},...,z_{n-1},z)\right\Vert
_{L^{2}(\left[ 0,T\right] ^{2(n-1)})}^{2}=0
\]%
$z$-a.e. Hence $\mathbb{P}_{1}-$a.e.  
\[
Y(t,x)=\mathbb{E}_{\mathbb{P}_{2}}\left[ Y(t,x)\right] +I_{1}(f_{1})
\]%
$\mathbb{P}_{2}-$a.e., that is $\mathbb{P}_{1}-$a.e. 

\begin{equation*}
Y(t,x)=\mathbb{E}_{\mathbb{P}_{2}}\left[ Y(t,x)\right] +\int_{0}^{t}%
\int_{0}^{x}f\left( (t-u)(x-v)\left( \frac{1}{N}A-I_{N}\right) \right)
B(du,dv)\text{ }\mathbb{P}_{2}-\text{a.e.}  
\end{equation*}%
Similarly, we find that 
\[
\mathbb{E}_{\mathbb{P}_{2}}\left[ Y(t,x)\right] =f\left( tx\left( \frac{1}{N}%
A-I_{N}\right) \right) Y(0)\text{ }\mathbb{P}_{1}-\text{a.e.}
\]%
We also note that%
\[
\mathbb{E}_{\mathbb{P}}\left[ Y(t,x)\left\vert \mathcal{G}\right. \right] =%
\mathbb{E}_{\mathbb{P}_{2}}\left[ Y(t,x)\right] \text{ }\mathbb{P}-\text{a.e.%
}
\]%
Therefore, $\mathbb{P}-\text{a.e.}$
\[
Y(t,x)=f\left( tx\left( \frac{1}{N}A-I_{N}\right) \right)
Y(0)+\int_{0}^{t}\int_{0}^{x}f\left( (t-u)(x-v)\left( \frac{1}{N}%
A-I_{N}\right) \right) B(du,dv)\text{ }
\]%
Define $\displaystyle\left\Vert A\right\Vert =\sum_{j=1}^{N}a_{j}\text{.}$
Require that $Y_{i}(0)=y$ for all $i=1,...,N$, the sequence $a_{j},j\geq 1$
is uniformly bounded and that%
\[
\frac{1}{N}\left\Vert A\right\Vert =\frac{1}{N}\sum_{j=1}^{N}a_{j}\geq q>0%
\text{ }\mathbb{P}-\text{a.e.}
\]%
for all $N$ and some constant $q$ and 
\[
\frac{1}{N}\left\Vert A\right\Vert =\frac{1}{N}\sum_{j=1}^{N}a_{j}\underset{%
N\longrightarrow \infty }{\longrightarrow }a>0\text{ }\mathbb{P}-\text{a.e.}
\]%
On the other hand, we observe for $n\geq 0$ that%
\begin{align*}
&\left( \frac{1}{N}A-I_{N}\right) ^{n} \nonumber\\&=\sum_{j=0}^{n}\binom{n}{j}\left( 
\frac{1}{N}A\right) ^{j}(-I_{N})^{n-j}=(-1)^{n}I_{N}+\sum_{j=1}^{n}\binom{n}{%
j}\left( \frac{1}{N}A\right) ^{j}(-I_{N})^{n-j}\nonumber \\
&=(-1)^{n}I_{N}+\sum_{j=1}^{n}\binom{n}{j}\left\Vert \frac{1}{N}%
A\right\Vert ^{j-1}\frac{1}{N}A(-I_{N})^{n-j}\nonumber \\
&=(-1)^{n}I_{N}+\left( -(-1)^{n}+(-1)^{n}+\sum_{j=1}^{n}\binom{n}{j}%
\left\Vert \frac{1}{N}A\right\Vert ^{j}(-1)^{n-j}\right) \frac{1}{\left\Vert 
\frac{1}{N}A\right\Vert }\frac{1}{N}A\nonumber \\
&=(-1)^{n}\left(I_{N}-\frac{1}{\left\Vert A\right\Vert }A\right)+\left(\frac{1}{N}%
\left\Vert A\right\Vert -1\right)^{n}\frac{1}{\left\Vert A\right\Vert }A\text{.}
\end{align*}%
So,%
\begin{eqnarray*}
&&f\left( (t-u)(x-v)\left(\frac{1}{N}A-I_{N}\right)\right) =\sum_{n\geq 0}\frac{1}{%
(n!)^{2}}((t-u)(x-v)\left(\frac{1}{N}A-I_{N})\right)^{n} \\
&=&\sum_{n\geq 0}\frac{1}{(n!)^{2}}((t-u)(x-v))^{n}\left\{ (-1)^{n}\left(I_{N}-%
\frac{1}{\left\Vert A\right\Vert }A\right)+\left(\frac{1}{N}\left\Vert A\right\Vert
-1\right)^{n}\frac{1}{\left\Vert A\right\Vert }A\right\}  \\
&=&f(-(t-u)(x-v))\left(I_{N}-\frac{1}{\left\Vert A\right\Vert }A\right)+f((t-u)(x-v)\left(
\frac{1}{N}\left\Vert A\right\Vert -1)\right)\frac{1}{\left\Vert A\right\Vert }A%
\text{.}
\end{eqnarray*}%
Hence%
\begin{align*}
Y(t,x) &=f(-tx)\left(I_{N}-\frac{1}{\left\Vert A\right\Vert }A\right)y+f(tx)\left( 
\frac{1}{N}\left\Vert A\right\Vert -1\right) \frac{1}{\left\Vert
A\right\Vert }Ay \nonumber\\
&+\int_{0}^{t}\int_{0}^{x}\left\{ f(-(t-u)(x-v))\left( I_{N}-\frac{1}{%
\left\Vert A\right\Vert }A\right) \right. \nonumber \\
&\left. +f((t-u)(x-v)\left( \frac{1}{N}\left\Vert A\right\Vert -1\right) 
\frac{1}{\left\Vert A\right\Vert }A\right\} B(du,dv)\text{.}
\end{align*}%
The latter entails that%
\begin{eqnarray*}
X^{i,N}(t,x) &=&f(-tx)\left( 1-\frac{\left\Vert A\right\Vert }{\left\Vert
A\right\Vert }\right) y+f(tx)\left( \frac{1}{N}\left\Vert A\right\Vert
-1\right) y \\
&&+\int_{0}^{t}\int_{0}^{x}f(-(t-u)(x-v))dB_{i}(du,dv)+I_{i,N}\text{,}
\end{eqnarray*}%
where 
\begin{eqnarray*}
I_{i,N}:= &&-\sum_{j=1}^{N}\int_{0}^{t}\int_{0}^{x}f(-(t-u)(x-v))\frac{a_{j}%
}{\left\Vert A\right\Vert }dB_{j}(du,dv) \\
&+&\sum_{j=1}^{N}\int_{0}^{t}\int_{0}^{x}f(-(t-u)(x-v))\left( \frac{1}{N}%
\left\Vert A\right\Vert -1\right) \frac{a_{j}}{\left\Vert A\right\Vert }%
dB_{j}(du,dv),
\end{eqnarray*}%
for all $1\leq i\leq N$. Using the It\^{o} isometry, our assumptions on $%
a_{j},j\geq 1$ and dominated convergence, we find that%
\begin{eqnarray*}
\mathbb{E}\left[ \left\vert I_{i,N}\right\vert ^{2}\right]  &\leq &C\frac{N}{%
\left\Vert A\right\Vert ^{2}}\int_{0}^{t}\int_{0}^{x}(f(-(t-u)(x-v)))^{2}\left(
\frac{1}{N}\left\Vert A\right\Vert -2\right)^{2}dudv \\
&&\underset{N\longrightarrow \infty }{\longrightarrow }0\text{,}
\end{eqnarray*}%
where $C$ is a constant. Thus, for $N\longrightarrow \infty $ 
\[
Y^{i,N}(t,x)\overset{d}{\longrightarrow }Y(t,x):=f(tx)(a-1)y+\int_{0}^{t}%
\int_{0}^{x}f(-(t-u)(x-v))dB^{\ast }(du,dv),
\]%
for all $i\geq 1$, where $B^{\ast }$ is another Brownian sheet independent
of $W$, and where $Y(t,x)$, $0\leq t$, $x\leq T$ solves (by the same
reasoning as above) the conditional mean-field hyperbolic SPDE 
\[
Y(t,x)=y+\int_{0}^{x}\int_{0}^{t}\Big(a\E\left[ Y(s,u)\left\vert \mathcal{G}
\right. \right] -Y(s,u)\Big)dsdu+B^{\ast }(t,x).
\]



\section{Conditional McKean-Vlasov SPDE}

In this section, we will discuss the existence and uniqueness of the solution for a conditional McKean-Vlasov SPDE driven by a Brownian sheet. The following tools will be in force throughout the rest of the paper.

We will consider the conditional probability distribution \( \mu_{t,x} = \mathcal{L}(Y(t,x) | \mathcal{F}_{t,x}^{(1)}) \) of the solution \( Y(t,x) \) to equation (\ref{(2)}).  We will specify the space where these probability distributions are defined and introduce the weighted  norm specifically tailored for the space of measures.

We assume that $m \geq 2$ and we fix one of the Brownian sheets, say $B_1=B_1(t,x,\omega)$, with filtration $\{\mathcal{F}_{t,x}^{(1)}\}_{t,x\geq 0}$. We define $\mu_{t,x}=\mu_{t,x}(dy,\omega)$ to be regular conditional distribution of $Y(t,x)$ given $\mathcal{F}_{t,x}^{(1)}$. This means that $\mu_{t,x}(dy,\omega) $ is a Borel probability measure on $\RR^n$ for all $(t,x) \in [0,T]\times [0,X],\omega \in \Omega$ and
\begin{equation} \label{cond}
\int_{\mathbb{R}^n} g(y) \mu_{t,x}(dy,\omega)= \EE[g(Y(t,x)) | \mathcal{F}_{t,x}^{(1)}](\omega),
\end{equation}
 for all functions $g$ such that $\EE[ |g(Y(t,x)) |] < \infty$. We refer to Theorem 9 in Plotter \cite{Plotter}.

We consider spaces containing these conditional probability distributions:

\begin{definition}
Let $n$ be a given natural number. Then let $\mathbb{M}=\mathbb{M}^{n}$ be the pre-Hilbert space of random measures $\mu$ on
$\mathbb{R}^{n}$ equipped with the norm  
\begin{align*}\label{norm}
\left\Vert \mu\right\Vert _{\mathbb{M}}^{2}:= \mathbb{E}\left[\int_{\mathbb{R}^{n}}
|\hat{\mu}(y)|^{2}e^{-y^{2}}dy \right]\text{,}
\end{align*}
where $y=(y_1,y_2, ... ,y_n)\in \mathbb{R}^{n}$ and $\hat{\mu}$ is the Fourier transform of the measure $\mu$, i.e.%
\[%
\begin{array}
[c]{lll}%
\hat{\mu}(y) & := & {\int_{\mathbb{R}^{n}}
}e^{-ixy}\mu(dx);\quad y\in\mathbb{R}^{n},
\end{array}
\]
where $xy =x \cdot y = x_1 y_1 + x_2 y_2 + ... + x_n y_n$ is the scalar product in $\mathbb{R}^{n}$.

If $\mu,\eta\in
\mathbb{M}$, we define the inner product $\left\langle \mu
,\eta\right\rangle _{\mathbb{M}}$ by
\small
\[
\left\langle \mu,\eta\right\rangle _{\mathbb{M}}=\E\Big[
\int_{{\mathbb{R}^{n}}}\operatorname{Re}(\overline{\hat{\mu}}(y)\hat{\eta
}(y))e^{-y^{2}}dy\Big],
\] 
\normalsize
where $\operatorname{Re}(z)$ denotes the real part and $\bar{z}$ denotes the complex conjugate of the complex number $z$. 
\end{definition}
The space $\mathbb{M}$ equipped with the inner
product $\left\langle \mu,\eta\right\rangle _{\mathbb{M}}$
is a pre-Hilbert space. 
Moreover, we have the following estimate (Lemma 2.1 in Agram et al. \cite{APO}):
Let $Y_{1}$ and $Y_{2}$ be two $n$-dimensional random variables in
$L^{2}(\mathbb{P})$ with associated conditional probability distributions $\mu_{1}$ and $\mu_{2}$, given the information $\mathcal{F}_{t,x}^{(1)}$, respectively. Then we have
\begin{align}\label{EST}
\left\Vert \mu_{1}-\mu_{2}\right\Vert _{\mathbb{M}
}^{2} \leq  \pi\ \mathbb{E}[(Y_{1}-Y_{2})^{2}]\text{.}
\end{align}
This shows that it provides a bound on the Wasserstein distance. Specifically, for the second Wasserstein distance $(p=2)$, we have:
\begin{align*}
\left\Vert \mu_{1}-\mu_{2}\right\Vert _{\mathbb{M}
}^{2} \leq  \pi W_2 ^2(\mu_{1},\mu_{2})\text{.}
\end{align*}
Using either space of measures equipped with the norm or the Wasserstein distance leads to consistent insights about the differences between probability distributions, though each metric may offer different perspectives or be more suitable for specific applications.


Specifically, we assume that the process $Y=(Y_1, ... ,Y_n)^{T} \in \R^{n\times 1}$ satisfies the following equation:
\small
\begin{align}\label{CMV}
Y(z)&
=
Y(0) 
+\int_{R_z} \alpha(\zeta) d\zeta
+\int_{R_z}\beta(\zeta)B(d\zeta),
\end{align}
where $z=(t,x)$ and 
$$\alpha(z)=\alpha(z,Y(z),\mu_{z})=(\alpha_1, ... ,\alpha_n)^{T} \in \R^{n \times 1},$$ 
and  
$$\beta(z)=\beta(z,Y(z),\mu_{z})= (\beta_{k,\ell}(z))_{1\leq k \leq n,1 \leq \ell \leq m} \in \R^{n \times m}.$$\\
Here $\mu_{t,x}$ represents the regular conditional probability distribution of $Y(t,x)$.
To study the well-posedness of the conditional McKean-Vlasov SPDE  \eqref{CMV} driven by Brownian sheet, we impose the following set of assumptions on the coefficients $\alpha$ and $\beta$:
\begin{itemize}
    \item [(a)] $\alpha_k(z,y,\mu):[0,T]^2\times \mathbb{R}^n \times \mathbb{M}\rightarrow \mathbb{R}, \beta_k(z,y,\mu):[0,T]^2\times \mathbb{R}^{n} \times \mathbb{M}\rightarrow \mathbb{R}^{1 \times m}$ are $\mathcal{F}_{t,x}$-adapted processes.
      \item  [(b)] There exists a constant $C$, which may differ from line to line, such that for all $z\in \mathbb{R}^{2}_+,y,y'\in \mathbb{R}^n,\mu,\mu'\in \mathbb{M}$, we have
\small
\begin{align*}
       \left\vert \alpha_k(z,y,\mu)-\alpha(z,y',\mu')\right\vert^2 &+\left\vert \beta_k(z,y,\mu)-\beta(z,y',\mu')\right\vert ^2\\
       &\leq C (\left\vert y-y'  \right\vert^2  +\left\vert \left\vert \mu-\mu'  \right\vert\right\vert^2 _{\mathbb{M}}).
\end{align*} 
Additionally, we suppose 
\begin{align*} 
\left\vert \alpha_k(\zeta,y,\mu)\right\vert^2+\left\vert \beta_k(\zeta,y,\mu)\right\vert^2 \leq C(1+\left\vert y\right\vert^2+\left\vert \left\vert \mu\right\vert\right\vert_{\mathbb{M}}^2).
\end{align*} 
\end{itemize}
These assumptions are based on \cite{Yeh85}, which established the existence of weak solutions to SDEs in the plane with continuous coefficients, and  \cite{Yeh87}, which proved the uniqueness of strong solutions even with deterministic boundary processes.

Let $J_0$ be the Bessel function of order zero and $r_0\approx 1.4458$ be the first nonnegative zero of $J_0$, in the following sense:
$$
r_0=inf \left\{t>0:J_0(2\sqrt{t})=\sum_{j=0}^\infty \frac{(-1)^j}{{j!}^2} t^j=0 \right\}.
$$
We shall recall the two-parameter version of Gronwall's
Lemma in \cite{ZN}:

\begin{lemma}[Two-parameter Gronwall's Lemma]
   Let $f$ be a non-negative and bounded function. Suppose there exists a constant $C_0>0$ which satisfies $C_0|z|\leq r_0$, such that
   $$f(z) \leq C_0 \int_{{R}_z} f(\zeta) d\zeta. $$
  Then $f$ vanishes on $R_z$.
\end{lemma}
\begin{theorem}[Existence and uniqueness]
Under the above assumptions (a)-(b), the conditional McKean-Vlasov SPDE \eqref{CMV} has a unique strong solution. 
\end{theorem}

\dproof 
The proof is based on the Picard iteration argument as in the proof of the Propagation of Chaos for Space time Ornstein-Uhlenbeck SPDE in Section 3.\\
\textbf{Step 1. Uniqueness.} Suppose that we have two solutions $Y,Y'$ and set $\tilde{Y}=Y-Y'$, such that $\tilde{Y}$ satisfies
\begin{align*}
\tilde{Y}(t,x)&=\int_{{R}_z}\{\alpha(\zeta,Y(\zeta),\mu_{\zeta})-\alpha(\zeta,Y'(\zeta),\mu'_{\zeta})\}d\zeta\nonumber\\
&+\int_{{R}_z}\{\beta(\zeta,Y(\zeta),\mu_{\zeta})-\beta(\zeta,Y'(\zeta),\mu'_{\zeta})\}dB(\zeta).
\end{align*}
Taking the mean square yields
\begin{align*}
\mathbb{E}|\tilde{Y}(t,x)|^2&=\mathbb{E}|\int_{{R}_z}\{\alpha(\zeta,Y(\zeta),\mu_{\zeta})-\alpha(\zeta,Y'(\zeta),\mu'_{\zeta})\}d\zeta\nonumber\\
&+\int_{{R}_z}\{\beta(\zeta,Y(\zeta),\mu_{\zeta})-\beta(\zeta,Y'(\zeta),\mu'_{\zeta})\}dB(\zeta)|^2.
\end{align*}
Application of the triangle inequality together with the linearity of the expectation, gives
\begin{align*}
\mathbb{E}|\tilde{Y}(t,x)|^2&\leq 2\mathbb{E}|\int_{{R}_z}\{\alpha(\zeta,Y(\zeta),\mu_{\zeta})-\alpha(\zeta,Y'(\zeta),\mu'_{\zeta})\}d\zeta|^2\nonumber\\
&+2\mathbb{E}|\int_{{R}_z}\{\beta(\zeta,Y(\zeta),\mu_{\zeta})-\beta(\zeta,Y'(\zeta),\mu'_{\zeta})\}dB(\zeta)|^2.
\end{align*}
We use the Cauchy-Schwarz inequality for the $d\zeta$-integral and the isometry for the $dB$-integral and we get 
\begin{align*}
\mathbb{E}|\tilde{Y}(t,x)|^2&\leq 2|z|^2 \mathbb{E}\int_{{R}_z}|\alpha(\zeta,Y(\zeta),\mu_{\zeta})-\alpha(\zeta,Y'(\zeta),\mu'_{\zeta})|^2d\zeta\nonumber\\
&+2\mathbb{E}\int_{{R}_z}|\beta(\zeta,Y(\zeta),\mu_{\zeta})-\beta(\zeta,Y'(\zeta),\mu'_{\zeta})|^2d\zeta,
\end{align*}
where we have used the notation $|z|=tx.$\\
Using Jensen's inequality combined with the Lipschitz condition, gives
\begin{align*}
&\mathbb{E}|\tilde{Y}(t,x)|^2\nonumber\\&\leq |z|^2 (C+C\pi)^2\int_{{R}_z} \mathbb{E}|Y(\zeta)-Y'(\zeta)|^2d\zeta\nonumber+(C+C\pi)^2\int_{{R}_z} \mathbb{E}|Y(\zeta)-Y'(\zeta)|^2d\zeta\nonumber\\
&=(C+C\pi)^2(|z|^2+1)\int_{{R}_z} \mathbb{E}|\tilde{Y}(\zeta)|^2d\zeta.
\end{align*}
By Gronwall's Lemma, we get
$\mathbb{E}|\tilde{Y}(t,x)|^2=0.$\\
\textbf{Step 2. Existence.} Define $Y^0(z)=y$ and $Y^n(z)$ inductively with corresponding probability distributions $\delta_y$ and $\mu^{n}_{\zeta}=P_{Y^n(\zeta)}$ respectively, as follows
\begin{align*}
Y^{n+1}(z)&=y+\int_{{R}_z}\alpha(\zeta,Y^{n}(\zeta),\mu^{n}_{\zeta})d\zeta+\int_{{R}_z}\beta(\zeta,Y^{n}(\zeta),\mu^{n}_{\zeta})dB(\zeta).
\end{align*}
Similar computations as in the uniqueness case, for some constant $K>0$ depending on the Lipschitz constant, lead to
\begin{align*}
\mathbb{E}|Y^{n+1}(z)-Y^n(z)|^2&\leq K^2|z|^2\int_{{R}_z} \mathbb{E}|Y^{n}(\zeta)-Y^{n-1}(\zeta)|^2d\zeta.
\end{align*}
Repeating this procedure $n$-times, we get 
\small
\begin{align*}
&\mathbb{E}|Y^{n+1}(z)-Y^n(z)|^2\\
&\leq K^{2n}|z|^2\int_{{R}_{z}}\int_{{R}_{z_{n}}} \cdot \cdot \cdot (\int_{{R}_{z_{n}}}\mathbb{E}|Y^{1}(\zeta_{n+1})-Y^{0}(\zeta_{n+1})|^2d\zeta_{n+1})d\zeta_{{n}}
  \cdot \cdot \cdot d\zeta_{z_1}\\
 &\leq K^{2n}|z|^{2n} sup_{u \in {R}_{z}} \mathbb{E}|Y^{1}(u)|^2 x_n.
  \end{align*}
  Taking the sum, we have
  \small
\begin{align*}
\sum_{n=0}^
{\infty} \mathbb{E}|Y^{n+1}(z)-Y^n(z)|^2&\leq sup_{u \in {R}_{z}} \mathbb{E}|Y^{1}(u)|^2 \sum_{n=0}^
{\infty} (K|z|)^{2n} x_n<\infty,
\end{align*}
with $x_n=-\sum_{j=1}^
{n}\frac{(-1)^j}{(j!)^{2}}x_{n-j}$ and $K|z|<\sqrt{r_0}$. Thus $(Y^n)_n$ converges.
\fproof


\section{An integral stochastic Fokker-Planck  equation for the conditional law}

In this section we state and prove an integral equation for the conditional law of a time-space McKean-Vlasov equation.\\

In the following we let $D$ denote the derivative in the sense of distributions on the Banach space $\mathbb{M}$ of Radon measures $m(dy)$ on $\mathbb{R}^n$ equipped with the total variation norm (regarding $\mathbb{M}$ as a subspace of the space $\mathcal{S}'$ of tempered distributions on $\mathbb{R}^n$).

Let $\mathcal{F} = \{\mathcal{F}_z; z \in \mathbb{R}^2_+ \}$ be a two-parameter filtration. 
Recall that a 2-parameter process $\{M(z); z \in \mathbb{R}^2_+\}$ is a martingale with respect to the filtration $\mathcal{F}$ if

(i) for each $z \in \mathbb{R}^2_+$ , $M(z)$ is $\mathcal{F}_z$-measurable;

(ii) for each $z \in \mathbb{R}^2_+$, $M(z) \in L^1(\mathbb{P})$; and

(iii) whenever $s_1\leq t_1$ and $s_2\leq t_2$  are both in $\mathbb{R}^2_+$, $$\E[M(t_1,t_2)|\mathcal{F}_{(s_1,s_2)}] = M(s_1,s_2), \,\, \textnormal{a.s.}$$
Note that, for example, in \cite{WZ74, I2}  a generalization of the martingale property to multidimensional parameter spaces was considered.
In the following we let $\mathcal{F}^{(2)}(\zeta)$ denote the filtration generated by $(B_2(\xi),B_3(\xi), ... ,B_m(\xi))$; $\xi \leq \zeta,$ so that $$\mathcal{F}_{\zeta}=\mathcal{F}^{(1)}_{\zeta} \times \mathcal{F}^{(2)}_{\zeta}.$$
Recall that if $a=(a_{1},a_{2}),b=(b_{1},b_{2})$, then
$$a\vee b=(\max(a_{1},b_{1}),\max(a_{2},b_{2})).
$$
Moreover, 
\begin{equation*}
  \begin{aligned}
   I((a_1,a_2)\bar{\wedge} (b_1,b_2))=
\begin{cases}
        1 \quad & \text{if } a_1 \leq b_1\quad  \text{and } \quad a_2 \geq b_2, \\ 
        0 \quad & \text{otherwise. } 
    \end{cases}
  \end{aligned}
\end{equation*}

To obtain the stochastic Fokker-Planck integro-differential equation driven by a Brownian sheet for the conditional law. In the non-conditioning case, we will rely on the results from the  paper \cite{AOPT2}. Our approach will involve employing It\^o's formula for stochastic integrals in the plane and Fourier transforms of measures.
In \cite{AO}, the method of Fourier transform of measures was used to derive the stochastic Fokker-Planck equation for the conditional distribution of the McKean-Vlasov jump diffusion process.

The following theorem  extends the results of \cite{AO}  to the case of  SPDEs driven by time-space Brownian motion:

\begin{theorem}[Integral equation for the conditional law]
Let\\ $\mu_z=\mathcal{L}(Y(z) |\mathcal{F}_z^{(1)})$ be the conditional law of the process $Y(z)$ given by \eqref{Y2}. Then $\mu_z$ satisfies the following stochastic integro-partial differential equation:
\small
    \begin{align}
    &\mu_z - \mu_0\nonumber\\
    &=\int_{R_z}\Big\{ \sum_{k=1}^n - D_k [\alpha_k(\zeta)\mu_{\zeta} ]
+\tfrac{1}{2}\sum_{k,\ell =1}^n D_k D_{\ell} \big(\beta_k(\zeta) \beta_{\ell}^{T}(\zeta)\mu_{\zeta}\big) \Big\} d\zeta \nonumber\\
    &+\int_{R_z} \Big\{ \sum_{k=1}^n -D_k \big(\beta_{k,1}(\zeta) \mu_{\zeta}\big)\Big\}B_1(d\zeta)\nonumber\\
    &+\iint\limits_{R_{z}\times R_{z}} \Big\{ \sum_{k,\ell =1}^n D_k D_{\ell}\big(\beta_{k,1}(\zeta)\beta_{\ell,1}(\zeta^{\prime})\mu_{\zeta \vee \zeta'}\big)\Big\} B_1(d\zeta) B_1(d\zeta')\nonumber\\
    &+\iint\limits_{R_{z}\times R_{z}} \Big\{\sum_{k,\ell =1}^n D_k D_{\ell}[ \big(\beta_k(\zeta^{\prime}) \alpha_{\ell}(\zeta) + \beta_k(\zeta) \alpha_{\ell}(\zeta')\big)\mu_{\zeta \vee \zeta'}]\nonumber\\
&-\tfrac{1}{2}\sum_{k,\ell,p=1}^n  D_k D_{\ell} D_p 
 [\big(\beta_{\ell}(\zeta)\beta_p^T (\zeta)\beta_{k,1}(\zeta^{\prime}) + \beta_{\ell}(\zeta')\beta_p^T (\zeta')\beta_{k,1}(\zeta)\big) \mu_{\zeta \vee \zeta'}]\Big\} d\zeta B_1(d\zeta')\nonumber\\
    &+\iint\limits_{R_{z}\times R_{z}} I(\zeta \bar{\wedge} \zeta^{\prime}) \Big\{\sum_{k,\ell =1}^n D_k D_{\ell} \big(\alpha_k(\zeta^{\prime}) \alpha_{\ell}(\zeta)\mu_{\zeta \vee \zeta'}\big)\nonumber\\
&-\tfrac{1}{2}\sum_{k,\ell,p=1}^n D_k D_{\ell} D_p
\big[( \alpha_k(\zeta^{\prime}) \beta_{\ell}(\zeta)\beta_p^T (\zeta)+ \alpha_k(\zeta) \beta_{\ell}(\zeta^{\prime}) \beta_p^{T}(\zeta^{\prime}))\mu_{\zeta \vee \zeta'}\big] \nonumber\\
& +\tfrac{1}{4}  \sum_{k,\ell,p,q=1}^n D_k D_{\ell} D_p D_q \big(\beta_k(\zeta^{\prime})\beta_\ell^{T}(\zeta^{\prime}%
)\beta_p(\zeta)\beta_q^{T}(\zeta) \mu_{\zeta \vee \zeta'}\big)\Big\} d\zeta d\zeta'.\label{SPDE}
   \end{align}
\end{theorem}
\dproof 
For given $z=(t,x)$ let $\mathcal{F}_z^{(1)}$ be the sigma-algebra generated by $\{B_1(\zeta);\zeta \leq z \}$, where $\zeta =(\zeta_1,\zeta_2) \leq z$ means $\zeta_1 \leq t,\zeta_2 \leq x$.
Choose   $\psi \in C^2(\mathbb{R})^n $ with bounded derivatives, and with values in the complex plane $\mathbb{C}$.  
Then by the It\^o formula we have
\small
\begin{align*}
&\E[\psi(Y(z))|\mathcal{F}_z^{(1)}]\nonumber\\
&=\psi(Y(0))+\int_{R_{z}} \E\Big[\sum_{k=1}^n \frac{\partial f}{\partial y_k}(Y(\zeta))\alpha_k(\zeta
)|\mathcal{F}_z^{(1)}\Big]d\zeta\nonumber\\
&+\E \Big[\sum_{k=1}^n \frac{\partial f}{\partial y_k}(Y(\zeta))\beta_{k,1}(\zeta)|\mathcal{F}_z^{(1)}\Big]B_1(d\zeta)\\
&+\tfrac{1}{2}\int_{R_{z}}\E\Big[\sum_{k,\ell =1}^n \frac{\partial^2 f}{\partial y_k \partial y_{\ell}} (Y(\zeta))\beta_k(\zeta) \beta_{\ell}^{T}(\zeta)|\mathcal{F}_z^{(1)}\Big]d\zeta\\
& +\iint\limits_{R_{z}\times R_{z}}\E\Big[\sum_{k,\ell =1}^n \frac{\partial ^2 f}{\partial y_k \partial y_{\ell}}(Y(\zeta\vee\zeta^{\prime
}))\beta_{k,1}(\zeta)\beta_{\ell,1}(\zeta^{\prime})|\mathcal{F}_z^{(1)}\Big] B_1(d\zeta)B_1(d\zeta^{\prime})\\ 
&+\iint\limits_{R_{z}\times R_{z}}\E\Big[\sum_{k,\ell =1}^n \frac{\partial ^2 f}{\partial y_k \partial y_{\ell}}(Y(\zeta\vee\zeta^{\prime
})){\beta_{k,1}(\zeta^{\prime}) \alpha
_{\ell}(\zeta)}\\
&+\tfrac{1}{2}\sum_{k,\ell,p=1}^n \frac{\partial ^{(3)} f}{\partial y_k \partial y_{\ell} \partial y_p}(Y(\zeta\vee\zeta^{\prime}))
 \beta_{\ell}(\zeta)\beta_p^T (\zeta)\beta_{k,1}(\zeta^{\prime})|\mathcal{F}_z^{(1)}\Big]d\zeta B_1(d\zeta^{\prime})\\
&+\iint\limits_{R_{z}\times R_{z}}\E\Big[\sum_{k,\ell =1}^n \frac{\partial ^2 f}{\partial y_k \partial y_{\ell}}(Y(\zeta\vee\zeta^{\prime
})){\beta_{k,1}(\zeta)\alpha_{\ell}(\zeta^{\prime})}\\
&+\tfrac{1}{2}\sum_{k,\ell,p=1}^n \frac{\partial ^{(3)} f}{\partial y_k \partial y_{\ell} \partial y_p}(Y(\zeta\vee\zeta^{\prime}))
 \beta_{\ell}(\zeta^{\prime})\beta_p^T (\zeta^{\prime})\beta_{k,1}(\zeta)|\mathcal{F}_z^{(1)}\Big]B_1(d\zeta)d\zeta^{\prime}\\
& +\iint\limits_{R_{z}\times R_{z}}I(\zeta \bar{\wedge} \zeta^{\prime})\E\Big[ \Big\{\sum_{k,\ell =1}^n \frac{\partial ^2 f}{\partial y_k \partial y_{\ell}}(Y(\zeta\vee\zeta^{\prime
}))\alpha_k(\zeta^{\prime}) \alpha_{\ell}(\zeta)\\
&
+\tfrac{1}{2}\sum_{k,\ell,p=1}^n \frac{\partial ^{(3)} f}{\partial y_k \partial y_{\ell} \partial y_p}(Y(\zeta\vee\zeta^{\prime}))
\Big[\alpha_k(\zeta^{\prime}) \beta_{\ell}(\zeta)\beta_p^T (\zeta)+ \alpha_k(\zeta) \beta_{\ell}(\zeta^{\prime}) \beta_p^{T}(\zeta^{\prime})\Big] \\
& +\tfrac{1}{4} \sum_{k,\ell,p,q=1}^n \frac{\partial ^4 f}{\partial y_k \partial y_{\ell} \partial y_p \partial y_q}(Y(\zeta\vee\zeta^{\prime}))\beta_k(\zeta^{\prime})\beta_\ell^{T}(\zeta^{\prime}%
)\beta_p(\zeta)\beta_q^{T}(\zeta)\Big\}|\mathcal{F}_z^{(1)}\Big]d\zeta d\zeta^{\prime}.
\end{align*}
This can be written
\small
\begin{align}
    \E[\psi(Y(z))|\mathcal{F}^{(1)}_z]-\psi(y)
    &=\int_{R_z}\E\Big[ A_1\psi(Y(\zeta))|\mathcal{F}^{(1)}_z\Big]d\zeta+\int_{R_z} \E\Big[A_2\psi(Y(\zeta))|\mathcal{F}^{(1)}_z\Big]B_1(d\zeta)\nonumber\\
    &+\iint\limits_{R_{z}\times R_{z}}\E\Big[A_3(\psi(Y(\zeta \vee \zeta')))|\mathcal{F}_z^{(1)}\Big] B_1(d\zeta) B_1(d\zeta')\nonumber\\
    &+\iint\limits_{R_{z}\times R_{z}}\E\Big[A_4(\psi(Y(\zeta \vee \zeta')))|\mathcal{F}_z^{(1)}\Big] d\zeta B_1(d\zeta')\nonumber\\
    &+\iint\limits_{R_{z}\times R_{z}}\E\Big[A_5(\psi(Y(\zeta \vee \zeta')))|\mathcal{F}_z^{(1)}\Big] d\zeta d\zeta',\label{ex.13}
   \end{align}
  where
  \small
    \begin{align*}
    &A_1\psi(Y(\zeta))= \sum_{k=1}^n \frac{\partial \psi}{\partial y_k}(Y(\zeta))\alpha_k(\zeta
)+\tfrac{1}{2}\sum_{k,\ell =1}^n \frac{\partial^2 \psi}{\partial y_k \partial y_{\ell}} (Y(\zeta))\beta_k(\zeta) \beta_{\ell}^{T}(\zeta);\\
    &A_2\psi(Y(\zeta))= \sum_{k=1}^n \frac{\partial \psi}{\partial y_k}(Y(\zeta))\beta_{k,1}(\zeta);\\
    &A_3 \psi(Y(\zeta \vee \zeta'))=\sum_{k,\ell =1}^n \frac{\partial ^2 \psi}{\partial y_k \partial y_{\ell}}(Y(\zeta\vee\zeta^{\prime
}))\beta_{k,1}(\zeta)\beta_{\ell,1}(\zeta^{\prime});\\
    &A_4 \psi(Y(\zeta \vee \zeta'))=\sum_{k,\ell =1}^n \frac{\partial ^2 \psi}{\partial y_k \partial y_{\ell}}(Y(\zeta\vee\zeta^{\prime
})){\big[\beta_{k,1}(\zeta^{\prime}) \alpha_{\ell}(\zeta) + \beta_{k,1}(\zeta) \alpha_{\ell}(\zeta')}\big]\nonumber\\
&+\tfrac{1}{2}\sum_{k,\ell,p=1}^n \frac{\partial ^{(3)} \psi}{\partial y_k \partial y_{\ell} \partial y_p}(Y(\zeta\vee\zeta^{\prime}))
 \big[\beta_{\ell}(\zeta)\beta_p^T (\zeta)\beta_{k,1}(\zeta^{\prime}) + \beta_{\ell}(\zeta')\beta_p^T (\zeta')\beta_{k,1}(\zeta)\big];\\  
&A_5 \psi(Y(\zeta \vee \zeta'))=I(\zeta \bar{\wedge} \zeta^{\prime}) \Big\{\sum_{k,\ell =1}^n \frac{\partial ^2 \psi}{\partial y_k \partial y_{\ell}}(Y(\zeta\vee\zeta^{\prime
}))\alpha_k(\zeta^{\prime}) \alpha_{\ell}(\zeta)\nonumber\\
&+\tfrac{1}{2} \sum_{k,\ell,p=1}^n \frac{\partial ^{(3)} \psi}{\partial y_k \partial y_{\ell} \partial y_p}(Y(\zeta\vee\zeta^{\prime}))
\Big[\alpha_k(\zeta^{\prime}) \beta_{\ell}(\zeta)\beta_p^T (\zeta)+\alpha_k(\zeta) \beta_{\ell}(\zeta^{\prime}) \beta_p^{T}(\zeta^{\prime})\Big] \nonumber\\
& +\tfrac{1}{4} \sum_{k,\ell,p,q=1}^n \frac{\partial ^4 \psi}{\partial y_k \partial y_{\ell} \partial y_p \partial y_q}(Y(\zeta\vee\zeta^{\prime}))\beta_k(\zeta^{\prime})\beta_\ell^{T}(\zeta^{\prime}%
)\beta_p(\zeta)\beta_q^{T}(\zeta)\Big\}.
    \end{align*}
We can write \eqref{ex.13} as follows:
    \small
\begin{align}
    \E[\psi(Y(z))|\mathcal{F}^{(1)}_z]-\psi(y)
    &=\E\Big[\int_{R_z}\E\Big[ A_1\psi(Y(\zeta))|\mathcal{F}^{(1)}_{\zeta}\Big]d\zeta+\int_{R_z} \E\Big[A_2\psi(Y(\zeta))|\mathcal{F}^{(1)}_{\zeta}\Big]B_1(d\zeta)\nonumber\\
    &+\iint\limits_{R_{z}\times R_{z}}\E\Big[A_3(\psi(Y(\zeta \vee \zeta')))|\mathcal{F}_{\zeta \vee \zeta'}^{(1)}\Big] B_1(d\zeta) B_1(d\zeta')\nonumber\\
    &+\iint\limits_{R_{z}\times R_{z}}\E\Big[A_4(\psi(Y(\zeta \vee \zeta')))|\mathcal{F}_{\zeta \vee \zeta'}^{(1)}\Big] d\zeta B_1(d\zeta')\nonumber\\
    &+\iint\limits_{R_{z}\times R_{z}}\E\Big[A_5(\psi(Y(\zeta \vee \zeta')))|\mathcal{F}_{\zeta \vee \zeta'}^{(1)}\Big] d\zeta d\zeta' |\mathcal{F}_z^{(1)}\Big].\label{5.13}
   \end{align}       
For given $w \in \mathbb{R}^n$ we now apply this to the function
$$\psi(y)= \exp( - i y w);\quad y \in \mathbb{R}^n,$$
where $i=\sqrt{-1}, y=(y_1, ..., y_n), w=(w_1, ... ,w_n) \text{ and } yw=y_1w_1 + y_2w_2+ ... +y_n w_n.$\\
Then we get 
\small
 \begin{align}
  & A_1\psi(Y(\zeta))=  \Big\{ \sum_{k=1}^n (-i) w_k \alpha_k(\zeta
)+\tfrac{1}{2}\sum_{k,\ell =1}^n w_k w_{\ell} \beta_k(\zeta) \beta_{\ell}^{T}(\zeta) \Big\} \exp(-iwY(\zeta) ),\label{5.20}\\
    &A_2\psi(Y(\zeta))=  \Big\{ \sum_{k=1}^n (-i)w_k \beta_{k,1}(\zeta)\Big\} \exp(-iwY(\zeta)),\label{5.21}\\
  & A_3 \psi(Y(\zeta \vee \zeta'))= \Big\{ \sum_{k,\ell =1}^n w_k w_{\ell} \beta_{k,1}(\zeta)\beta_{k,1}(\zeta^{\prime}))\Big\} \exp(-iwY(\zeta\vee\zeta^{\prime
})),\label{5.22}
\end{align}
 \begin{align}
    &A_4 \psi(Y(\zeta \vee \zeta'))=\Big\{\sum_{k,\ell =1}^n w_k w_{\ell} {\big(\beta_{k,1}(\zeta^{\prime}) \alpha_{\ell}(\zeta) + \beta_{k,1}(\zeta) \alpha_{\ell}(\zeta')\big)},\nonumber\\
&+\tfrac{1}{2}\sum_{k,\ell,p=1}^n (-i) w_k w_{\ell} w_p 
 \big(\beta_{\ell}(\zeta)\beta_p^T (\zeta)\beta_{k,1}(\zeta^{\prime}) + \beta_{\ell}(\zeta')\beta_p^T (\zeta')\beta_{k,1}(\zeta)\big)\Big\} \exp(-iwY(\zeta\vee\zeta^{\prime
})),\label{5.23}
  \end{align}
 \begin{align}
&A_5 \psi(Y(\zeta \vee \zeta'))=I(\zeta \bar{\wedge} \zeta^{\prime}) \Big\{\sum_{k,\ell =1}^n w_k w_{\ell}\alpha_k(\zeta^{\prime}) \alpha_{\ell}(\zeta)\nonumber\\
&+\sum_{k,\ell,p=1}^n (-i)w_k w_{\ell} w_p
\Big[\tfrac{1}{2}\alpha_k(\zeta^{\prime}) \beta_{\ell}(\zeta)\beta_p^T (\zeta)+\tfrac{1}{2} \alpha_k(\zeta) \beta_{\ell}(\zeta^{\prime}) \beta_p^{T}(\zeta^{\prime})\Big] \nonumber\\
& +\tfrac{1}{4} \sum_{k,\ell,p,q=1}^n w_k w_{\ell} w_p w_q \beta_k(\zeta^{\prime})\beta_\ell^{T}(\zeta^{\prime}%
)\beta_p(\zeta)\beta_q^{T}(\zeta)\Big\} \exp(-iwY(\zeta\vee\zeta^{\prime
})).\label{5.24}
    \end{align}
In general we have
\begin{align*}
   & \E[g(Y(z)) e^{-iY(z)w}|\mathcal{F}_z^{(1)}]=\int_{\mathbb{R}^n} g(y) e^{-iyw}\mu_z(dy)=F[g(\cdot) \mu_z(\cdot)](w), \nonumber\\
   &\text{ and } \\
   &\E[h(Y(z \vee z'))e^{-iY(z \vee z')w}|\mathcal{F}_{z \vee z'}^{(1)}]= \int_{\mathbb{R}^n} h(\eta) e^{-i\eta w}\mu_{z \vee z'}(d\eta)=F[h(\cdot) \mu_{z \vee z'}(\cdot)](w),
\end{align*}
where $F[\cdot]$ denotes the Fourier transform. Put
\small
 \begin{align*}
    &a_1=a_1(w)= \Big\{ \sum_{k=1}^n (-i) w_k \alpha_k(\zeta
)+\tfrac{1}{2}\sum_{k,\ell =1}^n w_k w_{\ell} \beta_k(\zeta) \beta_{\ell}^{T}(\zeta) \Big\}   \nonumber\\
    &a_2=a_2(w)= \Big\{ \sum_{k=1}^n (-i)w_k \beta_{k,1}(\zeta)\Big\}  \nonumber\\
    &a_3=a_3(w)=\Big\{ \sum_{k,\ell =1}^n w_k w_{\ell} \beta_{k,1}(\zeta)\beta_{\ell,1}(\zeta^{\prime})\Big\}  \nonumber\\
    &a_4=a_4(w)=\Big\{\sum_{k,\ell =1}^n w_k w_{\ell} {\big(\beta_{k,1}(\zeta^{\prime}) \alpha_{\ell}(\zeta) + \beta_{k,1}(\zeta) \alpha_{\ell}(\zeta')\big)}\nonumber\\
&+\tfrac{1}{2}\sum_{k,\ell,p=1}^n (-i) w_k w_{\ell} w_p 
 \big(\beta_{\ell}(\zeta)\beta_p^T (\zeta)\beta_{k,1}(\zeta^{\prime}) + \beta_{\ell}(\zeta')\beta_p^T (\zeta')\beta_{k,1}(\zeta)\big)\Big\}  \nonumber\\  
&a_5=a_5(w)=I(\zeta \bar{\wedge} \zeta^{\prime}) \Big\{\sum_{k,\ell =1}^n w_k w_{\ell}\alpha_k(\zeta^{\prime}) \alpha_{\ell}(\zeta)\nonumber\\
&+\sum_{k,\ell,p=1}^n (-i)w_k w_{\ell} w_p
\Big[\tfrac{1}{2}\alpha_k(\zeta^{\prime}) \beta_{\ell}(\zeta)\beta_p^T (\zeta)+\tfrac{1}{2} \alpha_k(\zeta) \beta_{\ell}(\zeta^{\prime}) \beta_p^{T}(\zeta^{\prime})\Big] \nonumber\\
& +\tfrac{1}{4} \sum_{k,\ell,p,q=1}^n w_k w_{\ell} w_p w_q \beta_k(\zeta^{\prime})\beta_\ell^{T}(\zeta^{\prime}%
)\beta_p(\zeta)\beta_q^{T}(\zeta)\Big\}. 
\end{align*}
By combining this with \eqref{5.13} and  \eqref{5.20} -- \eqref{5.24} we get
\small
\begin{align*}
F[Y(z)](w) - F[\delta_{Y(0)}](w)
&=\E[\exp(-iwY(z))] - \exp(-iwY(0))\nonumber\\
    &=\E[\psi(Y(z))|\mathcal{F}^{(1)}_z]-\psi(y)\nonumber\\
    &=\int_{R_z}\E\Big[ a_1 |\mathcal{F}^{(1)}_{\zeta}\Big]\exp(-iwY(\zeta))d\zeta\nonumber\\
    &+\int_{R_z} \E\Big[a_2 |\mathcal{F}^{(1)}_{\zeta}\Big]\exp(-iwY(\zeta))B_1(d\zeta)\nonumber\\
    &+\iint\limits_{R_{z}\times R_{z}}\E\Big[a_3|\mathcal{F}_{\zeta \vee \zeta'}^{(1)}\Big]\exp(-iwY(\zeta \vee \zeta')) B_1(d\zeta) B_1(d\zeta')\nonumber\\
    &+\iint\limits_{R_{z}\times R_{z}}\E\Big[a_4|\mathcal{F}_{\zeta \vee \zeta'}^{(1)}\Big] \exp(-iwY(\zeta \vee \zeta'))d\zeta B_1(d\zeta')\nonumber\\
    &+\iint\limits_{R_{z}\times R_{z}}\E\Big[a_5|\mathcal{F}_{\zeta \vee \zeta'}^{(1)}\Big]\exp(-iwY(\zeta \vee \zeta')) d\zeta d\zeta', 
    \end{align*}
and this is equivalent to
    \small
    \begin{align*}
    F[Y(z)](w) - F[\delta_{Y(0)}](w)
    &=\int_{R_z}\int_{\mathbb{R}^n} a_1 \exp(-iwy)\mu_{\zeta}(dy) d\zeta\nonumber\\
    &+\int_{R_z} \int_{\mathbb{R}^n}a_2 \exp(-iwy)\mu_{\zeta}(dy)B_1(d\zeta)\nonumber\\
    &+\iint\limits_{R_{z}\times R_{z}}\int_{\mathbb{R}^n} a_3 \exp(-iwy)\mu_{\zeta \vee \zeta'}(dy) B_1(d\zeta) B_1(d\zeta')\nonumber\\
    &+\iint\limits_{R_{z}\times R_{z}}\int_{\mathbb{R}^n} a_4 \exp(-iwy)\mu_{\zeta \vee \zeta'}(dy) d\zeta B_1(d\zeta')\nonumber\\
    &+\iint\limits_{R_{z}\times R_{z}}\int_{\mathbb{R}^n} a_5 \exp(-iwy)\mu_{\zeta \vee \zeta'}(dy) d\zeta d\zeta'. 
   \end{align*} 
In other words,
   \small
 \begin{align*}
    &F[Y(z)](w) - F[\delta_{Y(0)}](w)\nonumber\\
    &=\int_{R_z}\Big\{ \sum_{k=1}^n (-i) w_k F[\alpha_k(\zeta)\mu_{\zeta} ](w)
+\tfrac{1}{2}\sum_{k,\ell =1}^n w_k w_{\ell} F[\beta_k(\zeta) \beta_{\ell}^{T}(\zeta)\mu_{\zeta}](w] \Big\} d\zeta \nonumber\\
    &+\int_{R_z} \Big\{ \sum_{k=1}^n (-i)w_k F[\beta_{k,1}(\zeta) \mu_{\zeta}](w)\Big\}B_1(d\zeta)\nonumber\\
    &+\iint\limits_{R_{z}\times R_{z}} \Big\{ \sum_{k,\ell =1}^n w_k w_{\ell} F[\beta_{k,1}(\zeta)\beta_{\ell,1}(\zeta^{\prime})\mu_{\zeta \vee \zeta'}](w)\Big\} B_1(d\zeta) B_1(d\zeta')\nonumber\\
    &+\iint\limits_{R_{z}\times R_{z}} \Big\{\sum_{k,\ell =1}^n w_k w_{\ell}F[{\big(\beta_{k,1}(\zeta^{\prime}) \alpha_{\ell}(\zeta) + \beta_{k,1}(\zeta) \alpha_{\ell}(\zeta')\big)}\mu_{\zeta \vee \zeta'}](w)\nonumber\\
&+\tfrac{1}{2}\sum_{k,\ell,p=1}^n (-i) w_k w_{\ell} w_p 
 F[\big(\beta_{\ell}(\zeta)\beta_p^T (\zeta)\alpha_{k}(\zeta^{\prime}) + \beta_{\ell}(\zeta')\beta_p^T (\zeta')\alpha_{k}(\zeta)\big) \mu_{\zeta \vee \zeta'}](w)\Big\} d\zeta B_1(d\zeta')\nonumber\\
    &+\iint\limits_{R_{z}\times R_{z}} I(\zeta \bar{\wedge} \zeta^{\prime}) \Big\{\sum_{k,\ell =1}^n w_k w_{\ell} F[\alpha_k(\zeta^{\prime}) \alpha_{\ell}(\zeta)\mu_{\zeta \vee \zeta'}](w)\nonumber\\
&+\sum_{k,\ell,p=1}^n (-i)w_k w_{\ell} w_p
F\Big[\tfrac{1}{2} \alpha_k(\zeta^{\prime}) \beta_{\ell}(\zeta)\beta_p^T (\zeta)+\tfrac{1}{2} \alpha_k(\zeta) \beta_{\ell}(\zeta^{\prime}) \beta_p^{T}(\zeta^{\prime})\mu_{\zeta \vee \zeta'}\Big](w) \nonumber\\
& +\tfrac{1}{4} \sum_{k,\ell,p,q=1}^n w_k w_{\ell} w_p w_q F\Big[\beta_k(\zeta^{\prime})\beta_\ell^{T}(\zeta^{\prime}%
)\beta_p(\zeta)\beta_q^{T}(\zeta) \mu_{\zeta \vee \zeta'}\Big](w)\Big\} d\zeta d\zeta' .
   \end{align*}
   Recall that by the properties of the Fourier transforms, we have
\begin{align*}
    i w_k F[\alpha(y) \mu_{\zeta}(dy)](w)=F[D_k(\alpha(y) \mu_{\zeta}(dy))](w); \text{ where } D_k = \frac{\partial}{\partial y_k}
\end{align*}
and similar for higher order derivatives.
Therefore the above can be written
\small
\begin{align*}
    &F[\mu_z](w) - F[\mu_0](w)\nonumber\\
    &=\int_{R_z}\Big\{F[ \sum_{k=1}^n - D_k [\alpha_k(\zeta)\mu_{\zeta} ]](w)
+\tfrac{1}{2}F[\sum_{k,\ell =1}^n D_k D_{\ell} [\beta_k(\zeta) \beta_{\ell}^{T}(\zeta)\mu_{\zeta}]](w) \Big\} d\zeta \nonumber\\
    &+\int_{R_z} \Big\{ F[\sum_{k=1}^n -D_k [\beta_{k,1}(\zeta) \mu_{\zeta}]](w)\Big\}B_1(d\zeta)\nonumber\\
    &+\iint\limits_{R_{z}\times R_{z}} \Big\{ F[\sum_{k,\ell =1}^n D_k D_{\ell}[\beta_{k,1}(\zeta)\beta_{\ell,1}(\zeta^{\prime})\mu_{\zeta \vee \zeta'}]](w)\Big\} B_1(d\zeta) B_1(d\zeta')\nonumber\\
    &+\iint\limits_{R_{z}\times R_{z}} \Big\{F[\sum_{k,\ell =1}^n D_k D_{\ell}[ {\big(\beta_{k,1}(\zeta^{\prime}) \alpha_{\ell}(\zeta) + \beta_{k,1}(\zeta) \alpha_{\ell}(\zeta')\big)}\mu_{\zeta \vee \zeta'}]](w)\nonumber\\
&+\tfrac{1}{2}F[\sum_{k,\ell,p=1}^n - D_k D_{\ell} D_p 
 [\big(\beta_{\ell}(\zeta)\beta_p^T (\zeta)\alpha_{k}(\zeta^{\prime}) + \beta_{\ell}(\zeta')\beta_p^T (\zeta')\alpha_{k}(\zeta)\big) \mu_{\zeta \vee \zeta'}]](w)\Big\} d\zeta B_1(d\zeta')\nonumber\\
    &+\iint\limits_{R_{z}\times R_{z}} \Big\{F[I(\zeta \bar{\wedge} \zeta^{\prime}) \sum_{k,\ell =1}^n D_k D_{\ell} \big(\alpha_k(\zeta^{\prime}) \alpha_{\ell}(\zeta)\mu_{\zeta \vee \zeta'}\big)](w)\nonumber\\
&+\tfrac{1}{2}F[\sum_{k,\ell,p=1}^n -D_k D_{\ell} D_p
\big(\alpha_k(\zeta^{\prime}) \beta_{\ell}(\zeta)\beta_p^T (\zeta)+ \alpha_k(\zeta) \beta_{\ell}(\zeta^{\prime}) \beta_p^{T}(\zeta^{\prime})\mu_{\zeta \vee \zeta'}\big)](w) \nonumber\\
& +\tfrac{1}{4} F[ \sum_{k,\ell,p,q=1}^n D_k D_{\ell} D_p D_q \big(\beta_k(\zeta^{\prime})\beta_\ell^{T}(\zeta^{\prime}%
)\beta_p(\zeta)\beta_q^{T}(\zeta) \mu_{\zeta \vee \zeta'}\big)](w)\Big\} d\zeta d\zeta' .
   \end{align*}
   Since the Fourier transform describes the measure uniquely, we obtain \eqref{SPDE}.
   \fproof

\section{A stochastic partial differential equation for the conditional law }   
Next we prove the corresponding SPDE for the conditional law. To this end, the following result is useful:

\begin{lemma}{(Lemma 5.1 in \cite{AOPT2})} \label{Lemma}
Suppose $F=F_p$ has the form
\begin{align*}
       F(\zeta,\zeta',y,\mu)= D^p[f(\zeta,y,\mu)g(\zeta',y,\mu)\mu_{\zeta \vee \zeta'}];\quad p=1,2,3,4.
\end{align*}
Then
\small
    \begin{align*}
        &\frac{\partial^2}{\partial t \partial x} \iint\limits_{R_{z}\times R_{z}} I(\zeta \bar{\wedge} \zeta')F(\zeta,\zeta',y,\mu)d\zeta d\zeta'\\ &=D^p\Big[\Big(\int_0^t \int_0^x f((\zeta_1,x),y,\mu)g((t,\zeta'_2),y,\mu))d\zeta_1 d\zeta'_2\Big)\mu_{t,x}\Big]. 
    \end{align*}
In particular, if $f(\zeta,y,\mu)=f(\zeta)$ and $g(\zeta,y,\mu)=g(\zeta)$do not depend on $y$ and $\mu$, we get
  \begin{align*}
        \frac{\partial^2}{\partial t \partial x} \iint\limits_{R_{z}\times R_{z}} I(\zeta \bar{\wedge} \zeta')F(\zeta,\zeta',y,\mu)d\zeta d\zeta'=\Big(\int_0^t \int_0^x f(\zeta_1,x)g(t,\zeta'_2)d\zeta_1 d\zeta'_2\Big)D^p\mu_{t,x}, 
    \end{align*}  
and if $f,g$ are constants, we get
\begin{align*}
        \frac{\partial^2}{\partial t \partial x} \iint\limits_{R_{z}\times R_{z}} I(\zeta \bar{\wedge} \zeta')F(\zeta,\zeta',y,\mu)d\zeta d\zeta'=txfg D^p\mu_{t,x}. 
\end{align*}
\end{lemma}

\begin{theorem} \label{SPDE2}
Let $Y(z)$ be given by \eqref{Y2} and let $\mu_z=\mathcal{L}(Y(z) |\mathcal{F}_z^{(1)})$ be the conditional law of $Y(z)$. Then $\mu_z=\mu_{t,x}(y)$ satisfies the following stochastic partial differential equation:

    \begin{align}\label{SPDE 7.1}
    \frac{\partial^2}{\partial t \partial x}\mu_{t,x}(y)=A^{*}\mu_{t,x}(y)\nonumber\\
    \end{align}
    where $A^{*}=A_y^{*}$ is the integro-differential operator acting on  $y$ given by
    \small
    \begin{align}
    A^{*}\mu_z:=&\sum_{k=1}^n - D_k [\alpha_k(z)\mu_{z} ]
+\tfrac{1}{2}\sum_{k,\ell =1}^n D_k D_{\ell} \big(\beta_k(z) \beta_{\ell}^{T}(z)\mu_{z}\big) \nonumber\\
    &+\Big\{ \sum_{k=1}^n -D_k \big(\beta_{k,1}(z) \mu_{z}\big)\Big\}\diamond \overset{\bullet}{B}_1(z)\nonumber\\
    &+\Big\{\int_{R_{z}} \Big( \sum_{k,\ell =1}^n D_k D_{\ell}\big(\beta_{k,1}(z)\beta_{\ell,1}(\zeta^{\prime})\mu_{z}\Big) B_1(d\zeta')\Big\}\diamond \overset{\bullet}{B}_1(z)\nonumber\\
    &+\Big\{\int_{R_{z}} \Big( \sum_{k,\ell =1}^n D_k D_{\ell}\big(\beta_{k,1}(\zeta)\beta_{\ell,1}(z)\mu_{z}\Big) B_1(d\zeta)\Big\}\diamond \overset{\bullet}{B}_1(z)\nonumber\\
    &+\int_{R_{z}} \Big\{\sum_{k,\ell =1}^n D_k D_{\ell}[ \big(\beta_k(\zeta^{\prime}) \alpha_{\ell}(z) + \beta_k(z) \alpha_{\ell}(\zeta')\big)\mu_{z}]\nonumber\\
&-\tfrac{1}{2}\sum_{k,\ell,p=1}^n  D_k D_{\ell} D_p 
 [\big(\beta_{\ell}(z)\beta_p^T (z)\beta_{k,1}(\zeta^{\prime}) + \beta_{\ell}(\zeta')\beta_p^T (\zeta')\beta_{k,1}(z)\big) \mu_{z}]\Big\}d\zeta' \diamond \overset{\bullet}{B}_1(z)\nonumber\\
 &+\int_{R_{z}} \Big\{\sum_{k,\ell =1}^n D_k D_{\ell}[ {\big(\beta_{k,1}(z) \alpha_{\ell}(\zeta) + \beta_{k,1}(\zeta) \alpha_{\ell}(z)\big)}\mu_{z}]\nonumber\\
&-\tfrac{1}{2}\sum_{k,\ell,p=1}^n  D_k D_{\ell} D_p 
 [\big(\beta_{\ell}(\zeta)\beta_p^T (\zeta)\beta_{k,1}(z) + \beta_{\ell}(z)\beta_p^T (z)\beta_{k,1}(\zeta)\big) \mu_{z}]\Big\} d\zeta \diamond \overset{\bullet}{B}_1(z)\nonumber\\
    &+ \sum_{k,\ell =1}^n D_k D_{\ell}\Big[ \int_0^t \int_0^x\big(\alpha_k(t,\zeta_2^{\prime}) \alpha_{\ell}(\zeta_1,x)d\zeta_1 d\zeta_2'\big)\mu_{t,x}\Big]\nonumber\\
&-\tfrac{1}{2} \sum_{k,\ell,p=1}^n D_k D_{\ell} D_p \Big[(\int_0^t \int_0^x 
\big(\alpha_k(t,\zeta_2^{\prime}) \beta_{\ell}(\zeta_1,x)\beta_p^T (\zeta_1,x)\nonumber\\&+ \alpha_k(\zeta_1,x) \beta_{\ell}(t,\zeta_2^{\prime}) \beta_p^{T}(t,\zeta_2^{\prime})\big)d\zeta_1 d\zeta_2' ) \mu_{t,x}\Big] \nonumber\\
& +\tfrac{1}{4}  \sum_{k,\ell,p,q=1}^n D_k D_{\ell} D_p D_q \Big[\int_0^t \int_0^x \big(\beta_k(t,\zeta_2^{\prime})\beta_\ell^{T}(t,\zeta_2^{\prime}%
)\beta_p(\zeta_1,x)\beta_q^{T}(\zeta_1,x)\big) d\zeta_1 d\zeta_2'  \mu_{t,x}\Big].\label{SPDE3}
   \end{align}    
\end{theorem}

\dproof
First, note that if $\zeta \in R_z$ then $\mu_{\zeta \vee z}=\mu_z.$ Next, note that in general we have
\begin{align*}
    \frac{\partial^2}{\partial t \partial x} \iint\limits_{R_{z}\times R_{z}}f(\zeta,\zeta') d\zeta d\zeta'=\int_{R_z} f(z,\zeta') d\zeta' +\int_{R_z} f(\zeta,z) d\zeta.
\end{align*}
Using this and Lemma \ref{Lemma} we get Theorem \ref{SPDE2} by differentiating \eqref{SPDE}.
We omit the details.
\fproof

\begin{corollary}
Suppose the coefficients $\alpha,\beta$ do not depend on $y$. Then     
 \small
    \begin{align}
    \frac{\partial^2}{\partial t \partial x}\mu_{t,x}(y)
    &=\sum_{k=1}^n - \alpha_k(z) D_k [\mu_{z} ]
+\tfrac{1}{2}\sum_{k,\ell =1}^n \beta_k(z) \beta_{\ell}^{T}(z)D_k D_{\ell} [\mu_{z}] \nonumber\\
    &+\Big\{ \sum_{k=1}^n -\beta_{k,1}(z)D_k[\mu_{z}]\Big\}\diamond \overset{\bullet}{B}_1(z)\nonumber\\
    &+\Big\{\int_{R_{z}} \Big( \sum_{k,\ell =1}^n \beta_{k,1}(z)\beta_{\ell,1}(\zeta^{\prime})D_k D_{\ell}[\mu_{z}]\Big) B_1(d\zeta')\Big\}\diamond \overset{\bullet}{B}_1(z)\nonumber\\
    &+\Big\{\int_{R_{z}} \Big( \sum_{k,\ell =1}^n \beta_{k,1}(\zeta)\beta_{\ell,1}(z) D_k D_{\ell}[\mu_{z}] B_1(d\zeta)\Big\}\diamond \overset{\bullet}{B}_1(z)\nonumber\\
    &+\int_{R_{z}} \Big\{\sum_{k,\ell =1}^n \big(\beta_k(\zeta^{\prime}) \alpha_{\ell}(z) + \beta_k(z) \alpha_{\ell}(\zeta')\big)D_k D_{\ell}[ \mu_{z}]\nonumber\\
&-\tfrac{1}{2}\sum_{k,\ell,p=1}^n \big(\beta_{\ell}(z)\beta_p^T (z)\beta_{k,1}(\zeta^{\prime}) + \beta_{\ell}(\zeta')\beta_p^T (\zeta')\beta_{k,1}(z)\big) D_k D_{\ell} D_p 
 [ \mu_{z}]\Big\}d\zeta' \diamond \overset{\bullet}{B}_1(z)\nonumber\\
 &+\int_{R_{z}} \Big\{\sum_{k,\ell =1}^n {\big(\beta_{k,1}(z) \alpha_{\ell}(\zeta) + \beta_{k,1}(\zeta) \alpha_{\ell}(z)\big)}D_k D_{\ell}[ \mu_{z}]\nonumber\\
&-\tfrac{1}{2}\sum_{k,\ell,p=1}^n  \big(\beta_{\ell}(\zeta)\beta_p^T (\zeta)\beta_{k,1}(z) + \beta_{\ell}(z)\beta_p^T (z)\beta_{k,1}(\zeta)\big)D_k D_{\ell} D_p 
 [ \mu_{z}]\Big\} d\zeta \diamond \overset{\bullet}{B}_1(z)\nonumber\\
    &+ \sum_{k,\ell =1}^n (\int_0^t \int_0^x\big(\alpha_k(t,\zeta_2^{\prime}) \alpha_{\ell}(\zeta_1,x)d\zeta_1 d\zeta_2') D_k D_{\ell} \mu_{t,x} \nonumber\\
&-\tfrac{1}{2} \sum_{k,\ell,p=1}^n (\int_0^t \int_0^x 
\big(\alpha_k(t,\zeta_2^{\prime}) \beta_{\ell}(\zeta_1,x)\beta_p^T (\zeta_1,x)+ \alpha_k(\zeta_1,x) \beta_{\ell}(t,\zeta_2^{\prime}) \beta_p^{T}(t,\zeta_2^{\prime})\big)d\zeta_1 d\zeta_2') D_k D_{\ell}  D_p \mu_{t,x}\nonumber\\
& +\tfrac{1}{4}  \sum_{k,\ell,p,q=1}^n (\int_0^t \int_0^x \big(\beta_k(t,\zeta_2^{\prime})\beta_\ell^{T}(t,\zeta_2^{\prime}%
)\beta_p(\zeta_1,x)\beta_q^{T}(\zeta_1,x)\big) d\zeta_1 d\zeta_2') D_k D_{\ell} D_p D_q  \mu_{t,x}.\label{SPDE3_1}
   \end{align}          
\end{corollary}

\begin{corollary}\label{m_col}
    Assume that $\mu_{t,x}(y)$ is absolutely continuous with respect to Lebesgue measure $dy$, with Radon-Nikodym derivative
    \begin{align}
     m_{t,x}(y)= \frac{\mu_{t,x}(dy)}{dy}.   
    \end{align}
    Then
    \begin{align*}
    \frac{\partial^2}{\partial t \partial x}m_{t,x}(y)=A^{*}m_{t,x}(y)\nonumber\\
    \end{align*}
\end{corollary}
\begin{corollary}
Suppose the coefficients $\alpha,\beta$ are constants. Then
\small
    \begin{align*}
    &\frac{\partial^2}{\partial t \partial x}\mu_{t,x}(y)\nonumber\\
    &=\sum_{k=1}^n - \alpha_k D_k [\mu_{z} ]
+\tfrac{1}{2}\sum_{k,\ell =1}^n \beta_k \beta_{\ell}^{T}D_k D_{\ell} [\mu_{z}] \nonumber+\Big\{ \sum_{k=1}^n -\beta_{k,1} D_k [\mu_{z}]\Big\}\diamond \overset{\bullet}{B}_1(z)\nonumber\\
    &+2 \Big( \sum_{k,\ell =1}^n \beta_{k,1}\beta_{\ell,1} D_k D_{\ell}[\mu_{z}]\Big)\diamond B_1(z)\diamond \overset{\bullet}{B}_1(z)\nonumber\\
 &+tx \Big\{4\sum_{k,\ell =1}^n {\beta_{k,1} \alpha_{\ell}}D_k D_{\ell}[ \mu_{z}]\nonumber-2\sum_{k,\ell,p=1}^n  \beta_{\ell}\beta_p^T \beta_{k,1} D_k D_{\ell} D_p 
 [ \mu_{z}]\Big\}  \diamond \overset{\bullet}{B}_1(z)\nonumber\\
    &+ tx \sum_{k,\ell =1}^n \alpha_k \alpha_{\ell} D_k D_{\ell} [\mu_{t,x}] \nonumber-tx  \sum_{k,\ell,p=1}^n  
\alpha_k \beta_{\ell}\beta_p^T D_k D_{\ell}  D_p [\mu_{t,x}]\nonumber\\
& +\tfrac{1}{4} tx  \sum_{k,\ell,p,q=1}^n \beta_k\beta_\ell^{T}%
\beta_p \beta_q^{T} D_k D_{\ell} D_p D_q  [\mu_{t,x}].
   \end{align*}      
\end{corollary}

\section{Application to partial observation control}
Suppose we have a \emph{signal process} $Y(z)$ described by an equation of the form \eqref{Y2}:
\small
\begin{align*}\label{Ycont}
\text{(signal process)  }&Y(z)=Y(0)\nonumber\\
&+\int_{R_z} 
\alpha(\zeta,Y(\zeta),\mu_{\zeta},u(\zeta))d\zeta
+\int_{R_z}
\beta(\zeta,Y(\zeta),\mu_{\zeta},u(\zeta))B(d\zeta),
\end{align*}
 but now with coefficients $\alpha=(\alpha_1, \alpha_2, ... ,\alpha_n)\,\, \textnormal{and}\, \beta=(\beta_{i,j})_{1\leq i \leq n, 1\leq j\leq m}$ depending in addition on a control process $u=u(\zeta).$  Suppose we want to find a control $u(\cdot)$ which maximizes a given performance functional, defined by
\begin{align}
    J(u)=\E[\int_0^T \int_0^X \ell(\zeta,Y(\zeta),u(\zeta)) d\zeta + k(Y(T,X))]
\end{align}
for given utility functions $\ell$, $k$ and given terminal $T>0$,  $X>0$.  We assume that we know the dynamics of the system, in the sense that we know the functions $\alpha$, $\beta$, $\ell$, $k$,
but we only have indirect information about the state $Y(z)$ through the \emph{observation process} $G(t,x)$ given by the following equation
\begin{align*}
\text{(observation process)}\hskip 1.3cm   dG(z)=  dB_1(z)
\end{align*}
given by the first component  $B_1(\zeta);\zeta \leq z$ of the $m$-dimensional Brownian sheet $B=(B_1, B_2, ... , B_m)$.

Hence our control process $u(z)$ is required to be adapted to the filtration $\mathcal{F}_z^{G}= \mathcal{F}_z^{(1)}$ generated by the process $B_1(z)$. We let $\mathcal{A}_1$ denote the set of such processes. \\
Our \emph{partial observation problem} is the following:
\begin{problem}[Partial observation control problem]
    Find $u^{*} \in \mathcal{A}_1$ such that
 \begin{align}
     J(u^{*})= \sup_{u \in \mathcal{A}_1} J(u). \label{part1}
 \end{align}   
\end{problem}
We can rewrite this as a full observation problem as follows:\\
In general our best estimate of, say, $f(Y)$ at time-space $z=(t,x)$ is the conditional expectation 
$$\mathbb{E}[f(Y(z)) | \mathcal{F}_z^{(1)}]=\int_{\R^n} f(y) \mu_z(dy)$$
where $\mu_z(dy)$ is given by \eqref{SPDE3}.
Using this, we can write the performance as follows:
\begin{align*}
J(u)&=\E\Big[\int_0^T \int_0^X \ell(\zeta,Y(\zeta),u(\zeta))d\zeta + k(Y(T,X))\Big]\nonumber\\
&=\E\Big[\int_0^T \int_0^X \E[\ell(\zeta,Y(\zeta),u(\zeta))|\mathcal{F}_{\zeta}^{(1)}]d\zeta + \E[k(Y(T,X))|\mathcal{F}_{(T,X)}^{(1)}]\Big]\nonumber\\
&=\E\Big[\int_0^T \int_0^X \E[\ell(\zeta,Y(\zeta),v)|\mathcal{F}_{\zeta}^{(1)}]_{v=u(\zeta)}d\zeta + \E[k(Y(T,X))|\mathcal{F}_{(T,X)}^{(1)}]\Big]\nonumber\\
&=\E\Big[\int_0^T \int_0^X \int_{\R^n}\ell(\zeta,y,v) \mu_{\zeta}(dy)_{v=u(\zeta)}d\zeta + \int_{\R^n}k(y) \mu_{(T,X)}(dy)\Big]\nonumber\\
&=\E\Big[\int_0^T \int_0^X \int_{\R^n}\ell(\zeta,y,u(\zeta)) \mu_{\zeta}(dy)d\zeta + \int_{\R^n}k(y) \mu_{(T,X)}(dy)\Big]\nonumber\\
\end{align*}
We can summarize this as follows:
\begin{theorem}
The \emph{partial observation} problem \eqref{part1} is equivalent to the following \emph{full observation problem}:
\begin{problem}\label{full}
    Find $u^{*} \in \mathcal{A}_1$ such that
 \begin{align*}
     \widetilde{J}(u^{*})= \sup_{u \in \mathcal{A}_1} \widetilde{J}(u)
 \end{align*}
 where 
 \begin{align*}
  \widetilde{J}(u)=\E\Big[\int_0^T \int_0^X \int_{\R^n}\ell(\zeta,y,u(\zeta)) \mu_{\zeta}(dy)d\zeta + \int_{\R^n}k(y) \mu_{(T,X)}(dy)\Big]  
 \end{align*}
 and $\mu_{\zeta}(dy)$ is given by the SPDE \eqref{SPDE 7.1}.
\end{problem}
\end{theorem}
\begin{remark}
    Note that the SPDE \eqref{SPDE 7.1} is driven by $B_1(d\zeta)$ only. Therefore Problem \ref{full} is a full observation problem.
\end{remark}


\end{document}